\documentclass[12pt,a4paper]{article}

\usepackage{amsmath, amssymb, setspace, graphics}

\def\W{{\mathcal{W}}}
\def\bydef{\mathrel{:=}}
\def\w{{\mathbf{w}}}
\def\qed{{\hfill $\diamondsuit$}}
\def\CP{{{\mathbb C}{\rm P}}}
\def\Aut{{\rm Aut\,}}
\def\Emb{{\rm Emb\,}}
\def\hG{{\widehat G}}
\def\tG{{\widetilde G}}
\def\id{{\rm id}}
\def\Tr{{\rm Tr\,}}

\def\C{{\mathbb C}}

\def\ie{{\it i.e.}}
\def \lmod#1\rmod {\vphantom{#1}\left|\smash{#1}\right|}
\newcommand{\name}[1] {\operatorname{\rm #1}}
\newcommand{\pder}[2] {\frac{\partial #1}{\partial #2}}
\newcommand{\pdertwo}[3] {\frac{\partial^2 #1}{\partial #2 \partial #3}}

\usepackage{theorem}

\newtheorem{theorem}{Theorem}
\newtheorem{proposition}{Proposition}[section]
\newtheorem{corollary}[proposition]{Corollary}
\newtheorem{lemma}[proposition]{Lemma}
\newtheorem{conjecture}[proposition]{Conjecture}
{\theorembodyfont{\rmfamily}
\newtheorem{definition}[proposition]{Definition}
\newtheorem{example}[proposition]{Example}
\newtheorem{remark}[proposition]{Remark}
\newtheorem{notation}[proposition]{Notation}
}

\newcommand{\addtext}[1]{{\bfseries $>\!\!>$#1$<\!\!<$}\marginpar{V}}

\include{epsf}

\title{Cycle factorizations and one-faced graph embeddings}

\author{Yurii Burman\thanks{Independent University of Moscow, B.\
Vlassievsky per.~11, and Higher School of Economics, Moscow, Myasnitskaya,
20; e-mail: burman@mccme.ru. Supported by the CRDF grant RUM1-2895-MO-07,
INTAS grant 05-7805, RFBR grants 08-01-00110-a and NSh-709.2008.1, and the
HSE Scientific Foundation grant 08-01-0019.},
Dimitri Zvonkine\thanks{Poncelet laboratory, CNRS, Independent
University of Moscow, B.\  Vlassievsky per.~11; e-mail:
dimitri.zvonkine@gmail.com. Partially supported by the ANR grant
ANR-05-BLAN-0029-01
(Geometry and Integrability in Mathematical Physics).}}

\date{}

\begin{document}

\maketitle

 \begin{abstract}
Consider factorizations into transpositions of an $n$-cycle in the
symmetric group $S_n$. To every such factorization we
assign a monomial in variables $w_{ij}$ that retains the transpositions
used, but forgets their order. Summing over all possible factorizations of
$n$-cycles we obtain a polynomial that happens to admit a closed
expression. From this expression we deduce a formula for
the number of $1$-faced embeddings of a given graph.
 \end{abstract}

\section{Introduction and main results}
The Hurwitz problem is the problem of
counting ramified coverings of surfaces with prescribed ramification types.
It is a classical problem presently enjoying a regain of interest due to
its discovered relations with Young-Mills models, matrix integrals, and
intersection theory on moduli spaces~\cite{ELSV,GroTay,BouMar}.

\subsection{Hurwitz numbers and Hurwitz polynomials}

In this paper we not only count ramified coverings, but actually retain an
important part of the structure of each ramified covering. As a
consequence, the answer to the Hurwitz problem will be a polynomial in many
variables rather than a number.

We concentrate on a particular case of {\em pseudo-polynomial Morse}
coverings of the sphere by an arbitrary surface. In other words, we consider
degree~$n$ ramified coverings of the sphere by a genus~$g$ surface with
full ramification over one point labeled $\infty$ (\ie, $\infty$ has a
unique preimage) and simple ramifications over $n+2g-1$ other points (\ie,
each of them has one double and $n-2$ simple preimages). The problem of
counting pseudo-polynomial Morse coverings arises as a particular case
in~\cite{PouSch}, \cite{ShShVa}, and~\cite{ELSV} and is completely solved.
Our approach allows us, however, not only to count the coverings, but also
to obtain new information on the structure of the set of these coverings.
In particular, we draw some nontrivial consequences on the number of
one-faced graph embeddings. In the Appendix and
partly in Section~\ref{Sec:GrEmb} we also deal with coverings that are
Morse, but not necessarily pseudo-polynomial.

Choosing a base point on the sphere (different from the branch points) and
numbering the preimages of the base point, we obtain a description of every
ramified covering in terms of permutations. The monodromy of a
pseudo-polynomial Morse covering over~$\infty$ is an $n$-cycle in $S_n$,
while its monodromies over the other branch points are transpositions. The
monodromies determine the covering uniquely up to isomorphism. Therefore we
are actually interested in describing the set of factorizations of
$n$-cycles into $n+2g-1$ transpositions.

\begin{definition} \label{Df:CycleFact}
A list of transpositions $\tau_1, \dots, \tau_{n+2g-1} \in S_n$ such that
the product $\tau_{n+2g-1} \cdots \tau_1$ is an $n$-cycle is called a {\em
genus $g$ cycle factorizion}.
\end{definition}

Denote by $\C_n[\w]$ the ring of polynomials in (commuting) variables
$w_{ij}$, $1 \leq i,j \leq n$, $i \ne j$, modulo the relations $w_{ij} =
w_{ji}$. (We could have restricted ourselves to the variables $w_{ij}$ with
$i <j$, but sometimes it is convenient to use $w_{ij}$ without bothering to
know if $i$ is smaller or greater than~$j$.)

To every transposition $\tau \in S_n$ we assign the variable $w(\tau) =
w_{ij}$, where $i$ and $j$ are the elements permuted by~$\tau$.

 \begin{definition} \label{Df:HurwPoly}
The {\em Hurwitz polynomial} $P_{g,n}(\w)$ is defined by
$$
P_{g,n}(\w) = \sum w(\tau_1) \cdots w(\tau_{n+2g-1}),
$$
where the sum is taken over all cycle factorizations of genus~$g$.

The {\em Hurwitz number} $h_{g,n}$ is defined by
$$
h_{g,n} = \frac1{n!} P_{g,n}(\bf 1),
$$
where $P_{g,n}(\bf 1)$ is the result of the substition $w_{ij}=1$
for all $i,j$ in $P_{g,n}$.
 \end{definition}

The Hurwitz number is the number of all pseudo-polynomial Morse coverings,
each  covering being counted with weight $1/\mbox{(number of its
symmetries)}$.

 \begin{example}
For $n=2$ we have $P_{g,n} = w_{12}^{2g+1}$, $h_{g,n} = 1/2$.
 \end{example}

The Hurwitz polynomial $P_{g,n}$ is the main object of our study.
From each cycle factorization $P_{g,n}$ retains the transpositions
that compose it, but forgets their order.

The set of cycle factorizations is invariant under the action of $S_n$.
Hence $P_{g,n}$ is $S_n$-invariant under the renumberings of indices of the
varibles $w_{ij}$. Therefore it is a natural idea to combine similar terms
and represent their sum as a graph.

Recall that an automorphism of a graph is a permutation of its half-edges
preserving the relations ``to belong to the same edge'' and ``to have a
common vertex''.

 \begin{definition}
Let $G$ be a graph with no loops. We denote by $G_n(\w) \in \C_n[\w]$ the
polynomial in variables $w_{ij}$, obtained as follows. Label the vertices
of~$G$ with distinct numbers from $1$ to~$n$ in all possible ways (if $G$
has $v \leq n$ vertices there are $n!/(n-v)!$ labellings; if $v > n$ there
are no possible labellings). To each labelling assign the product of the
variables $w_{ij}$ corresponding to the edges. Sum the obtained monomials
over all labellings. Divide by the number $\lmod\Aut G\rmod$ of
automorphisms of~$G$.
 \end{definition}

 \begin{remark}
The same graph $G$ represents a polynomial $G_n\in C_n[\w]$ for each~$n$.
If $G$ has more than~$n$ vertices, then $G_n(\w)=0$.
 \end{remark}

 \begin{example} $\vphantom{a}$\par
For $G = \epsfbox[-1 5 5 25]{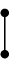}$ we have
$G_n(\w) = \frac12 \sum\limits_{i \not= j} w_{ij} =
\sum\limits_{i<j} w_{ij}$.

For $G  = \epsfbox[-1 5 10 25]{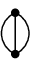}$ we have
$G_n(\w) = \frac1{12} \sum\limits_{i \not= j} w_{ij}^3 =
\frac16 \sum\limits_{i<j} w_{ij}^3$.

For $G  = \epsfbox[-1 10 7 37]{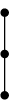}$ we have
$G_n(\w) = \sum\limits_{
\substack{i \not= j \not= k \\
i<k}
} w_{ij} w_{jk}.$
\end{example}

Let $\pi_n: \C_n[\w] \to \C_{n-1}[\w]$ be the substitution $w_{1n} = \dots
= w_{n-1,n} = 0$.

 \begin{definition}
The algebra $\W$ is the projective limit of spaces of $S_n$-invariants in
$\C_n[\w]$ with respect to the projections $\pi_n$.
 \end{definition}

A family of $S_n$-invariant polynomials $P_n(\w) \in \C_n[\w]$ is an
element of $\W$ if their degrees are uniformly bounded and
$$
\left.P_n\right|_{w_{1n} = \dots = w_{n-1,n}=0} = P_{n-1}
$$
for all~$n$. For example, for any graph $G$ the sequence $G_n(\w)$ defines
an element of $\W$, to be denoted as $G(\w)$ or just $G$ if it does not
lead to ambiguity. Moreover, it is clear that the elements $G(\w)$ for all
$G$ span $\W$.

 \begin{example}
We have
 %*
 \begin{align*}
\epsfbox[-1 5 6 20]{palochka.eps}^2 &= \left(\sum_{i<j} w_{ij} \right)^2
= \sum_{i<j} w_{ij}^2 + 2\sum_{\substack{i \ne j \ne k \\ i<k}}
w_{ij} w_{jk} +
2 \!\!\!\! \sum_{\substack{i<j, i<k<l \\j \ne l, j \ne k}}
\!\!\!\!w_{ij} w_{kl}\\
&= 2 \Bigl(
\epsfbox[-1 5 10 25]{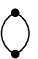} +
\epsfbox[-1 10 7 37]{2palochka.eps} +
\epsfbox[-1 5 12 25]{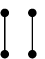}
\Bigr)
 \end{align*}
 %*
 \end{example}

Our aim is now to give an explicit expression for~$P_{g,n}(\w)$. To do this
we need to introduce some notation.

 \begin{definition} \label{Def:T}
We let $T_n(\w) = \sum_t t_n(\w)$, where the sum is taken over all
trees~$t$ with $n$ vertices.
 \end{definition}

 \begin{example}
We have
$$
\begin{array}{ccccl}
T_1(\w) &=& 1 && \mbox{(by convention)},\\
T_2(\w) &=& \epsfbox[-1 5 3 25]{palochka.eps} &=& w_{12},\\
T_3(\w) &=& \epsfbox[-1 10 3 37]{2palochka.eps} &=&
w_{12}w_{13}+w_{12}w_{23}+w_{13}w_{23},\\
T_4(\w) &=& \epsfbox[-2 18 45 55]{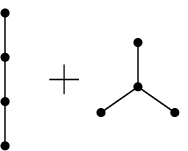} & = & \mbox{$16$ terms.}
\end{array}
$$
 \end{example}

\bigskip

\begin{notation}
We will repeatedly use the function
$$
\phi(t) = \frac{\sinh(t/2)}{t/2}
$$
and its logarithm
$$
\ln \phi = \sum_{g \geq 0} \frac{B_{2g}}{(2g)!\; (2g)} t^{2g},
$$
where $B_{2g}$ are the Bernoulli numbers.
\end{notation}

 \begin{notation} \label{Def:ArR}
Denote by $A_n$ the following $n \times n$ matrix:
$$
A_n = \left(
\begin{array}{ccccc}
\sum w_{1i} & -w_{12} & -w_{13} & \cdots & -w_{1n}\\
-w_{12} & \sum w_{2i} & -w_{23} & \cdots & -w_{2n} \\
-w_{13} & -w_{23} & \ddots & \ddots & \vdots\\
\vdots & \vdots & \ddots & \ddots & -w_{n-1,n}\\
-w_{1n} & -w_{2n} & \dots &-w_{n-1,n} & \sum w_{ni}
\end{array}
\right).
$$
 \end{notation}
 \begin{definition}
Introduce the power series
 %*
 \begin{equation}\label{Eq:Defr}
r_n(\w) = \Tr \ln \phi(A_n) = \sum_{g \geq 0} \frac{B_{2g}}{(2g)!\; (2g)}
\Tr A_n^{2g}
 \end{equation}
 %*
and
 %*
 \begin{equation}\label{Eq:DefR}
R_n(\w) = \exp (r_n(\w)) = \det \phi(A_n).
 \end{equation}
 %*
Both $r_n$ and $R_n$ are even. We denote by $r_{g,n}$ and $R_{g,n}$,
respectively, their homogeneous parts of degree~$2g$.
 \end{definition}

 \begin{proposition} \label{Prop:indepn}
For any given~$g$, $R_g = (R_{g,n})_{n \geq 1}$ and $r_g = (r_{g,n})_{n
\geq 1}$ define elements of~$\W$.
 \end{proposition}

\paragraph{Proof.} We have
$$
\left. A_n^{2g}\right|_{w_{1n} = \dots = w_{n-1,n} = 0} = \left(
 \begin{array}{c|c}
A_{n-1}^{2g} & 0 \\
\cline{1-2} 0 & 0
 \end{array}\right).
$$
The second expression for $r_n$ in \eqref{Eq:Defr} implies then
$\left.r_n\right|_{w_{1n} = \dots = w_{n-1,n} = 0} = r_{n-1}$. The result
about $R_n$ follows.\qed

\begin{example}
We have $r_0 = 0$, $R_0 = 1$, $r_1 = R_1 = \epsfbox[0 10 62 25]{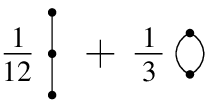}$,
 %*
 \begin{align*}
%&r_1 = R_1 = \epsfbox[0 10 62 25]{R1.eps},\\
&r_2 = \epsfbox[0 52 260 75]{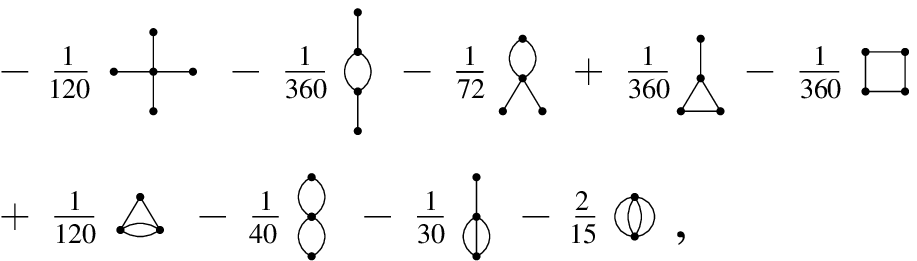}\\
&R_2 = \epsfbox[25 56 450 130]{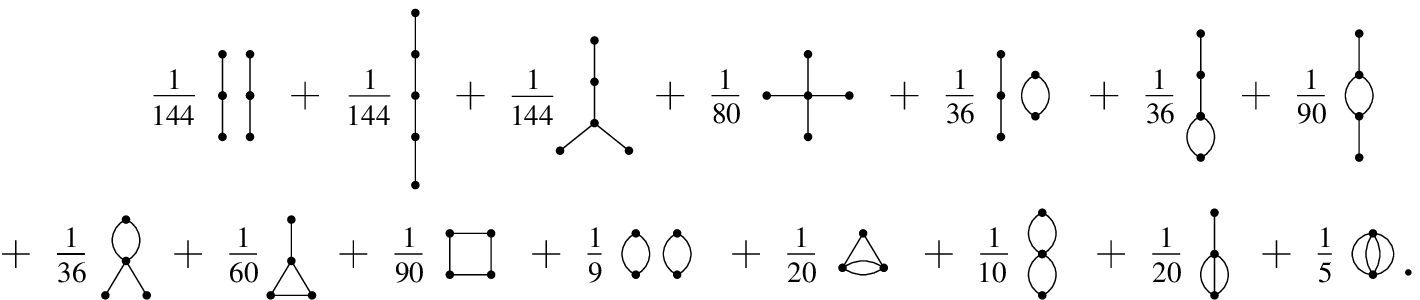}
 \end{align*}
 %*
\vspace{1.5cm}
\end{example}

 \begin{remark}
Let us write $r_g = \sum c_G G$ and $R_g = \sum C_G
G$, where $G$ are graphs. Then the sign of $c_G$ is equal to $(-1)^{b(G)}$
where $b(G)$ is the number of edges of $G$ that belong to at least one cycle
--- see the proof of Lemma \ref{Lem:union} for details.

 \begin{conjecture}
All the coefficients $C_G$ of the elements $R_g \in \W$ are positive.
 \end{conjecture}
 \end{remark}

The importance of $R_n$ is due to their role in the following theorem:

 \begin{theorem} \label{Thm:div}
We have
$$
P_{g,n}(\w) = (n+2g-1)! \; T_n(\w) \; R_{g,n}(\w).
$$
 \end{theorem}

 \begin{example}
We have $P_{0,n} = (n-1)! \; T_n$. The latter equality is interpreted as
follows: for every tree with~$n$ vertices, the product of the $n-1$
transpositions corresponding to its edges is an $n$-cycle {\em whatever the
ordering of the transpositions}. The coefficient $(n-1)!$ is simply the
number of possible orderings. \end{example}

\begin{corollary}{\rm\cite[Theorem 6]{ShShVa}} \label{Cor:Hurwitz}
The Hurwitz number $h_{g,n}$ equals
$$
h_{g,n} = \frac{(n+2g-1)!}{n!} \; n^{n-2} \; \rho_{g,n},
$$
where $\rho_{g,n}$ is the coefficient of $w^{2g}$ in the
power series $(\phi(nw))^{n-1}$.
\end{corollary}

\subsection{Graph embeddings}

An important application of our results is the counting of $1$-faced
embeddings of arbitrary graphs.

A connected graph is {\em embedded} into a closed oriented surface of
genus~$g$ if it is drawn on the surface with no intersections between the
edges (except at their endpoints) and if the edges cut the surface into
topological discs. These discs are called {\em faces}. Two embeddings of a
graph $G \to C_1$ and $G \to C_2$ are considered isomorphic if there exists
an orientation-preserving homeomorphism between $C_1$ and $C_2$ that makes
the triangular diagram commute:

\setlength\unitlength{1em}
\begin{center}
\
\begin{picture}(6,3.8)
\put(2.4,3.5){$G$}
\put(2.4,3.3){\vector(-2,-3){1.5}}
\put(3,3.3){\vector(2,-3){1.5}}
\put(0,0){$C_1$}
\put(4.4,0){$C_2$}
\put(1.2,.45){\vector(1,0){3}}
\end{picture}
\end{center}

Up to this isomorphism, giving an embedding of~$G$ is the same thing as
giving a {\em cyclic order} of half-edges issuing from every vertex of~$G$.
Indeed, if the embedding is given then the cyclic order is just the
counterclockwise order of the half-edges on the surface. Conversely, given
a cyclic order of half-edges at every vertex, it is possible to reconstruct
the faces of the covering and to glue a disc into each face thus obtaining
the surface of embedding.

If the valencies of the vertices of a connected graph $G$ are $k_1, \dots,
k_n$ then it has $\Emb G = \prod (k_i-1)!$ different embeddings.

For a detailed introduction into graph embeddings see~\cite{LanZvo}.

 \begin{definition}
Let $G$ be a connected graph with $2g$ independent cycles,
\ie, $\beta_1(G)=2g$.

A {\em decoration} of $G$ is the choice of several vertices of~$G$ and, for 
each of these vertices~$v$, the choice of a positive {\em even} number 
$k_v$ of half-edges adjacent to it, such that $\sum k_v = 2g$ and that if 
we erase the chosen half-edges the remaining part of~$G$ is contractible.

The {\em weight} of a decoration equals
$$
\frac1{2^{2g}} \prod_v \frac1{k_v+1}.
$$
 \end{definition}

 \begin{theorem} \label{Thm:spiders}
Let $G$ be a connected graph with $2g$ independent cycles.
The number of $1$-faced (that is, genus~$g$) embeddings of~$G$
divided by the total number of its embeddings is equal
to the sum of weights of all decorations of~$G$.
 \end{theorem}

\begin{example}
Consider the following graphs:
$$
G_1 = \epsfbox[-1 5 11 20]{tripleedge.eps},
\quad G_2 = \epsfbox[-1 1.5 23 20]{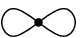},
\quad G_3 = \epsfbox[-1 1.5 29 20]{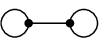}.
$$
They have, respectively, $(3-1)! \,(3-1)! = 4$, $(4-1)! = 6$,
and $(3-1)! \,(3-1)! = 4$ embeddings, out of which, respectively
$2$, $2$, and none are $1$-faced.

The decorations of $G_1$ and $G_2$ are shown in the figure.
Each of them has weight $1/12$. The graph $G_3$ has no decorations.

\begin{center}
\
\epsfbox{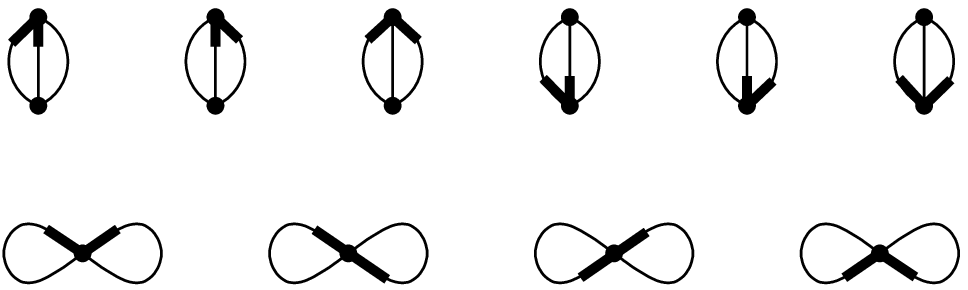}
\end{center}

For $G_1$ we have $6/12 = 2/4$, for $G_2$ we have $4/12 = 2/6$,
and for $G_3$ we have $0 = 0/4$.
\end{example}

\begin{example}
Let $G$ be the graph with $1$ vertex and $2g$ loops.
Then it has $2^{2g}$ decorations, each of weight $1/2^{2g}(2g+1)$.
It also has $(4g-1)!$ embeddings, of which $(4g-1)!/(2g+1)$
are $1$-faced: this is a particular case of
the Harer--Zagier formula~\cite{HarZag}. We have
$$
\frac{2^{2g}}{2^{2g}\, (2g+1)} = \frac{(4g-1)!/(2g+1)}{(4g-1)!}.
$$
\end{example}

\begin{example}
The graph in the figure has 96 decorations. One of them is shown
in the figure, while the others are obtained from it by
graph automorphisms. Every decoration has weight
$\frac1{2^6 \cdot 3 \cdot 5} = \frac1{960}$.

\begin{center}
\
\epsfbox{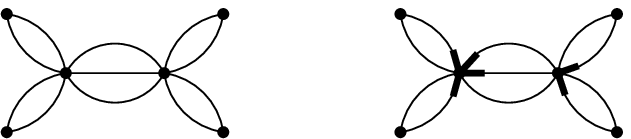}
\end{center}
Therefore this graph has $6!^2 \cdot \frac{96}{960} = 51840$ embeddings
into a genus~$3$ surface.
\end{example}

\section{Representations of $S_n$ and Frobenius formula}

The main goal of this section is to prove Theorem~\ref{Thm:div}.

Throughout the section we'll be using the following representations of the
symmetric group $S_n$

 \begin{itemize}
\item $\C^n$, where $S_n$ acts by permutations of coordinates (the
``permutation representation'').

\item $V \bydef \{(x_1, \dots, x_n) \mid \sum_i x_i = 0\}\subset \C^n$ ---
an $(n-1)$-dimensional irreducible representation (``geometric'' or
``defining'' representation of $S_n$ as a Coxeter group).

\item The exterior powers $\bigwedge^k V$, $0 \leq k \leq n-1$. These are 
also irreducible representations; $\bigwedge^k V$ 
corresponds to the ``hook'' Young diagram (see
\cite[page 48]{FuHa})
\begin{center}
\epsfbox{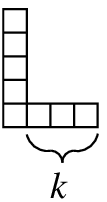}\;.
\end{center}

\item The group algebra $\C S_n$ where $S_n$ acts by left multiplication. A
classical theorem says that it is a direct sum of all the irreducible
representations of $S_n$ containing every representation $\lambda$ exactly
$\dim\lambda$ times.
 \end{itemize}

Let $C \in \C S_n$ be the sum of all $n$-cycles.
Let $B(\w) \in \C_n[\w] S_n$ be the sum
$$
B(\w) = \sum_{i <j} w_{ij} \tau_{ij},
$$
where $\tau_{ij}$ is the transposition interchanging $i$ and~$j$.

\begin{proposition} \label{prop:Frob}
$$
P_{g,n}(\w) = \frac1n \sum_{k=0}^{n-1} (-1)^k \; \Tr_{\bigwedge^k V}
\bigl( B(\w)^m \bigr),
$$
where $m= n+2g-1$.
\end{proposition}

\paragraph{Proof.} This is a particular case of the well-known Frobenius 
formula (see \cite{LanZvo}). For the reader's convenience we summarize the 
proof here.

For a permutation $\sigma \in S_n$, the trace $\Tr_{\C S_n} \sigma$ equals
$n!$ if $\sigma$ is the identity and $0$ otherwise. In other words,
$\frac{1}{n!}\Tr_{\C S_n} x$ is the coefficient of the identity permutation
in $x$. It follows that
$$
P_{g,n}(\w) = \frac1{n!} \Tr_{\C S_n} \bigl( B(\w)^m C \bigr)
= \frac1{n!} \sum_{\lambda} \dim \lambda \cdot \Tr_\lambda
\bigl( B(\w)^m C \bigr),
$$
where the sum is taken over the irreducible representations~$\lambda$.
Now, $C$ belongs to the center of $\C S_n$ and therefore, by Schur's
lemma, acts in every irreducible representation~$\lambda$
by a scalar. Using the ribbon rule \cite[Example 5 \S 3]{Macdonald}
for the character of an $n$-cycle, we see that this scalar
equals $(-1)^k (n-1)!/(\dim \lambda)$ for $\lambda = \bigwedge^k V$
and vanishes otherwise. Substituting this in the above
formula we obtain the equality of the proposition. \qed

 \begin{proposition} \label{prop:Liealg}
Let $\tau \in S_n$ be a transposition. Then the action of $\tau-1$ in
$\bigwedge^k \C^n  = \C^n \wedge \C^n \wedge \dots \wedge \C^n$ is given by
$$
\begin{array}{lcccccccc}
&\/\/ & (\tau - 1) &\wedge &\id &\wedge &\dots &\wedge& \id \\
+&& \id &\wedge &(\tau - 1)& \wedge& \dots& \wedge& \id \\
+ && &&&& \dots \\
+ &&\id& \wedge& \id&  \wedge& \dots& \wedge& (\tau - 1).
\end{array}
$$
The same formula is true for the exterior powers of the representation $V$.
 \end{proposition}

\paragraph{Proof.} Assume $\tau$ is the transposition~$(1,2)$.
Let $X = x_{i_1} \wedge \dots \wedge x_{i_p} \in \bigwedge^p \C^n$
be a wedge product that does not contain $x_1$ and $x_2$.
It is easy to check that both $\tau-1$ and the sum given
in the proposition act as follows (we have $p=k-2$, $k-1$, or $k$):
$$
\begin{array}{rcl}
X & \mapsto & X,\\
x_1 \wedge X & \mapsto & (x_2-x_1) \wedge X,\\
x_2 \wedge X & \mapsto & (x_1-x_2) \wedge X,\\
x_1 \wedge x_2 \wedge X & \mapsto & 2 \, x_2 \wedge x_1 \wedge X.
\end{array}
$$
The last statement follows from the fact that $V$ is a subrepresentation of
$\C^n$.\qed

 \begin{remark}
Proposition \ref{prop:Liealg} is crucial, although its proof is a trivial
check. One can formulate it by saying that $\tau-1$ is ``Lie-algebra-like''
for exterior powers of the representation~$\C^n$. Lie-algebra-like elements
form a vector subspace of $\C S_n$, and we believe they deserve a special
study. Examples include the difference of $3$-cycles $(123) - (132)$ and
two linear combinations of $4$-cycles: $(1234)-(1324)-(1423)+(1432)$ and
$(1243)-(1324)+(1342)-(1423)$.
 \end{remark}

Denote by $A(\w)$ the Lie-algebra-like element of $\C S_n$
$$
A(\w) = \sum w_{ij} (1-\tau_{ij}) = \sum w_{ij} - B(\w).
$$
The matrix of the action of $A(\w)$ in $\C^n$ is $A_n$ of Definition
\ref{Def:ArR}.

\begin{proposition} \label{Prop:eigenvalues}
Let $\sigma_1, \dots, \sigma_{n-1}$ be the eigenvalues of\/
$A(\w)$ in~$V$. Then the eigenvalues of
$B(\w)$ in $\bigwedge^k V$ are
$$
\frac{\pm \sigma_1 \pm \sigma_2 \pm \dots \pm \sigma_{n-1}}2
$$
with $k$ signs~``$-$'' and $n-1-k$ signs~``$+$''.
\end{proposition}

\paragraph{Proof.} Since $A(\w)$ is Lie-algebra-like, its eigenvalues in
$\bigwedge^k V$ are all sums $\sigma_{i_1} + \dots + \sigma_{i_k}$ for
distinct $i_1, \dots, i_k$. In particular, $A(\w)$ obviously acts by zero
in the trivial representation and by $\sum_{i=1}^{n-1} \sigma_i$ in the
($1$-dimensional) sign representation $\bigwedge^{n-1} V$. On the other
hand, in the sign representation every transposition acts by $-1$, so that
$A(\w)$ acts by $2 \sum_{i<j} w_{ij}$. Thus we have
$$
\sum_{i<j} w_{ij} = \frac12 \sum_{i=1}^{n-1} \sigma_i.
$$
This equality allows us to write the eigenvalue of $B(\w)$ in terms of
$\sigma_1, \dots, \sigma_{n-1}$ only, which gives the formula of the
proposition. \qed

 \begin{proposition} \label{Prop:kirchhof}
We have $\sigma_1 \cdots \sigma_{n-1} = nT_n(\w)$.
 \end{proposition}

\paragraph{Proof.} Let $v_1, \dots, v_n$ be the standard basis in $\C^n$; 
then the $(n-1)$-dimensional representation $V \subset \C^n$ is spanned by 
the vectors $u_2 \bydef v_2 - (v_1+\dots+v_n)/n, \dots, u_n \bydef v_n - 
(v_1+\dots+v_n)/n$. The matrix of the operator $A(\w)$ in this basis is 
then $(M_n)^{-1} \tilde A_n$ where 
$$
M_n = \frac{1}{n} \left(\begin{array}{cccc}
n & -1 & \dots & -1\\
-1 & n & \dots & -1\\
\hdotsfor4\\
-1 & -1 & \dots & n
\end{array}\right)
$$
and $\tilde A_n$ is obtained by erasing the first column and the first row 
from the matrix~$A_n$ of Notation~\ref{Def:ArR}. 

The Kirchhoff formula \cite{Kirchhoff} (a.k.a.\ the matrix-tree 
theorem) implies that $\det \tilde A_n = T_n(\w)$. Substituting $w_{ij} = 
1$ for all $i, j$ gives also $\det M_n = \frac{1}{n^{n-1}} T_n({\bf 1}) = 
1/n$.\qed

\bigskip

Now we are ready to prove two results announced in the introduction.

\paragraph{Proof of Theorem~\ref{Thm:div}.}
According to Proposition~\ref{Prop:eigenvalues}, we have
\begin{equation}\label{Eq:SumSign}
P_{g,n}(\w) = \frac{1}{n}\!\!\! 
\sum_{\varepsilon_1, \dots, \varepsilon_{n-1} \in \{-1,1\}}
\!\!\!\!\!\!\!\! \varepsilon_1 \dots \varepsilon_{n-1} \;
(\varepsilon_1 \sigma_1/2 + \dots + \varepsilon_{n-1}
\sigma_{n-1}/2)^{n+2g-1}.
 \end{equation}

This formula implies that $P_{g,n}$ is odd as a function of every
$\sigma_i$ (all the other $\sigma_j$ being fixed). Therefore, if we expand
the expression for $P_{g,n}$, all monomials $\sigma_1^{m_1} \dots
\sigma_{n-1}^{m_{n-1}}$ with at least one even $m_i$ cancel out. On the
other hand, if all the $m_i$ are odd, the contribution of every term
in~\eqref{Eq:SumSign} to the coefficient of this monomial is the same and
equal to
$$
\frac{(n+2g-1)!}{2^{n+2g-1}\, n \; m_1! \dots m_{n-1}!}.
$$
The total number of terms is $2^{n-1}$, so we finally get
$$
P_{g,n} \; =
\hspace{-1.5em}
\sum_{
\substack{\ \\
\ \\
m_1 + \cdots + m_{n-1} = n+2g-1\\
m_1, \dots, m_{n-1} \;\; {\rm odd}}
}
\hspace{-2em}
\frac{1}{2^{2g}\,n} \;
\frac{(n+2g-1)!}{m_1! \cdots m_{n-1}!} \;
\sigma_1^{m_1} \dots \sigma_{n-1}^{m_{n-1}},
\hspace{7em} \
$$
$$
\ \hspace{5em}
= (n+2g-1)! \;\; \frac{1}{n}\, \sigma_1 \cdots \sigma_{n-1}
\hspace{-1.5em}
\sum_{
\substack{\ \\
\ \\
p_1 + \cdots + p_{n-1} = g}
}
\hspace{-1.5em}
\frac{(\sigma_1/2)^{2p_1}}{(2p_1+1)!} \cdots
\frac{(\sigma_{n-1}/2)^{2p_{n-1}}}{(2p_{n-1}+1)!}
$$
By Proposition~\ref{Prop:kirchhof}, we have $\frac{1}{n}\, \sigma_1 \cdots 
\sigma_{n-1} = T_n(\w)$, so it remains to identify the last sum as 
$R_{g,n}$. Recall that
$$
\phi(\sigma) = \frac{\sinh(\sigma/2)}{\sigma/2}
=\sum_{p \geq 0} \frac{(\sigma/2)^{2p}}{(2p+1)!}.
$$
Thus we see that
$$
\sum_{
\substack{\ \\
\ \\
p_1, \dots, p_{n-1}}
}
\hspace{-.7em}
\frac{(\sigma_1/2)^{2p_1}}{(2p_1+1)!} \cdots
\frac{(\sigma_{n-1}/2)^{2p_{n-1}}}{(2p_{n-1}+1)!}
\; = \;
\prod_{i=1}^{n-1} \phi(\sigma_i)
\; = \;
\det\nolimits_V \phi(A_n).
$$
In the last determinant we can replace $V$ by $\C^n$. Indeed, the
difference between these two representations is the trivial representation,
where $A(\w)$ acts by~$0$. It remains to recall that $R_{g,n}$ was defined
as the degree $2g$ part of~$\det\nolimits_{\C^n} \phi(A_n)$. \qed

\paragraph{Proof of Corollary~\ref{Cor:Hurwitz}.} After the substitution
$w_{ij} = w$, the operator $A(\w)$ acts in~$V$ by the scalar~$nw$. In other
words, all the $\sigma_i$ are equal to $nw$, so we have
$$
\sum_{g \geq 0} \rho_{g,n} w^{2g}
\; = \; \det\nolimits_V \phi(A(\w))
= \; (\phi(nw))^{n-1}.
$$
\qed

\section{Graph embeddings}\label{Sec:GrEmb}

The aim of this section is to deduce several corollaries of
Theorem~\ref{Thm:div}, in particular to prove Theorem~\ref{Thm:spiders}. To
prove Theorem~\ref{Thm:spiders} we relate the coefficient of a graph in
$P_{g,n}$ with the number of $1$-faced embeddings of this graph and then
compute this coefficient using the right-hand side of the equality of
Theorem~\ref{Thm:div}.

First, explain the relation between the polynomials $P_{g,n}$ and graph
embeddings.

Recall that an embedding of a graph into a surface is defined, up to
equivalence, by the cyclic order of its half-edges at every vertex.

\begin{definition} \label{Def:numbedges}
Let $G$ be a connected graph without loops.
To a numbering of the edges of~$G$,
we assign the {\em embedding of the numbering} obtained
in the following way: the cyclic order
of the half-edges at each vertex is simply the increasing
order of the numbers of edges.
\end{definition}

\begin{remark}
In general, not every embedding can be obtained in this way: for instance,
a graph with $2$ vertices joined by $k$ edges has $(k-1)!^2$ embeddings and
only $k!$ numberings.
\end{remark}

\begin{remark} \label{Rem:preimage}
Consider a Morse covering of $\CP^1$ of degree~$n$ with~$s$
preimages of $\infty$. Let $\sigma \in S_n$ be its monodromy at
$\infty$ (a premutation with $s$ cycles). As we have seen in the introduction,
the covering determines (once we have chosen a set of
generators of $\pi_1(\C \setminus \{\mbox{branch points} \}$)
a factorization of $\sigma$ into transpositions. This factorization, in turn, 
determines a graph with numbered vertices and edges,
endowed with a natural embedding that we have just defined
(see Definition~\ref{Def:numbedges}).

There exists a more concise way to obtain the embedded
graph from a ramified covering. Let $z_1, \dots, z_{2g-2+n+s}$
be the branch points of the covering and $z_0$ a base point
(distinct from the branch points). On the plane draw a star
as follows:
\begin{center}
\
\epsfbox{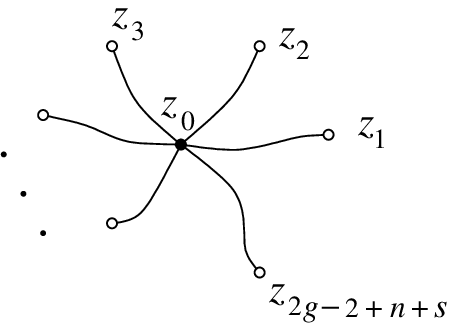}
\end{center}
This is equivalent to choosing a set of generators of 
$\pi_1(\C \setminus \{z_1, \dots, z_{2g-2+n+s} \})$.
The arrows of the star will be considered as 
half-edges. Now the procedure to obtain an embedded graph
from the ramified covering is simply the following:
take the preimage of the star and erase the dangling half-edges.
(Every $z_k$ has $n-2$ simple preimages and one double
preimage. Out of a simple preimage issues only one half-edge; 
these half-edges are called dangling. Out of a double
preimage issue two half-edges that form a complete edge.)

We leave it to the reader to check that the embedded
graph with numbred edges that we have obtained coincides
with the construction via factorizations into transpositions.
\end{remark}

Let $G$ be a connected graph without loops, with $n$ vertices and with
numbered edges. Consider the corresponding embedding and attribute a color
to each face. Follow the boundary of a face in the clockwise direction. For
each vertex of the boundary, if the edge preceding it has a larger number
than the edge following it, or the vertex is of valency one, then mark the
vertex with the color of the face.

 \begin{proposition}\label{Pp:CyclesEmb}
In the situation above, consider the product $\sigma \in S_n$ of the
transpositions corresponding to the edges in the order specified by the
numbering. Then the cycles of~$\sigma$ are in a one-to-one correspondence
with the faces of the embedding. The cycle corresponding to a face is
composed of the vertices bearing the color of the face.
 \end{proposition}

\paragraph{Proof.}
Introduce the following cyclic order on the vertices of a given color. Go
around the face of this color in the clockwise direction. As we meet two
successive edges whose numbers go in a nonincreasing order, the vertex
between them is declared the next with respect to the order.

Let us prove that (i)~each vertex of $G$ is colored exactly once, and
appears exactly once in the order corresponding to this color, (ii)~every
color is present, and (iii)~the vertices of a given color taken in the
corresponding order form a cycle of~$\sigma$.

(i)~The color of a vertex is simply the color of the face that
lies after its largest and before its smallest edge. It appears in the
order when we meet this pair of edges.

(ii)~If some color was absent it would mean that the
numbers of edges around the corresponding face increase
indefinitely.

(iii)~As we apply one by one the transpositions corresponding
to the edges, a vertex of a given color moves along the
face of the same color following a sequence of edges with increasing
numbers.
\begin{center}
\
\epsfbox{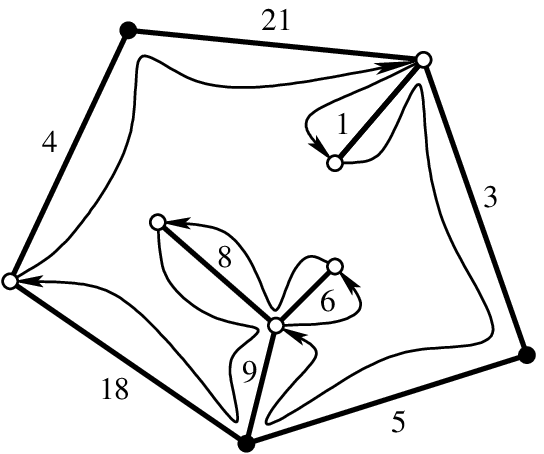}
\end{center}
In the end it arrives at the next vertex of the same color (with
respect to the order). The figure shows an example of a face, the
vertices of the corresponding color being represented in white.
\qed

\begin{remark}
Using the construction of Remark~\ref{Rem:preimage} we
can immediately conclude that the faces of the embedding
of~$G$ correspond to the cycles of~$\sigma$: indeed both
are in a natural one-to-one correspondence with the
poles of the covering.
\end{remark}

 \begin{remark}
In \cite{GouldYong} the relation between embeddings and numberings of edges
is studied in detail for the case when $G$ is a tree.
 \end{remark}

 \begin{corollary} \label{Cr:1Face}
Let $G$ be a connected graph without loops, with $n$ vertices, and with
$2g$ independent cycles, (\ie, its first Betti number equals $2g$). The
coefficient of $G$ in $P_{g,n}$ is  equal to the number of numberings of
the edges of~$G$ such that the corresponding embedding is $1$-faced.
 \end{corollary}

\paragraph{Proof.} Indeed, the Hurwitz polynomial enumerates
the factorizations of an $n$-cycle. \qed

 \begin{definition}
Let $G$ be a connected graph with at least 2 independent cycles. The
subgraph $\hG$ obtained from~$G$ by a repeated cutting off of valency~$1$
vertices is called the {\em skeleton} of~$G$. The vertices of~$G$ that have
valencies $\geq 3$ in~$\hG$ are called the {\em essential vertices} of~$G$.
 \end{definition}

\begin{center}
\
\epsfbox{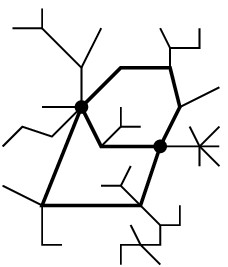}

The skeleton and the essential vertices of a graph.
\end{center}

\begin{definition}
A connected graph~$G$ with no loops and at least $2$ independent cycles
is called {\em long} if its essential vertices are never
connected by an edge.
\end{definition}

 \begin{proposition} \label{Prop:embeddings}
Let~$G$ be a long graph with $n$~vertices and $2g$~independent cycles. The
coefficient of~$G$ in $P_{g,n}/(n+2g-1)!$ equals the number of $1$-faced
embeddings of~$G$ divided by the number of all its embeddings.
 \end{proposition}

\paragraph{Proof.} Denote by $\Emb G$, $\Emb \hG$,
$\Emb\!_1 G$, and $\Emb\!_1 \hG$ the number of all embeddings
and of $1$-faced embeddings of $G$ and $\hG$.

The number of faces in the embedding of a long graph $G$ depends only on
the cyclic order $\alpha$ of the half-edges of the skeleton adjacent to the
essential vertices. Therefore we have
$$
\frac{\Emb\!_1 G}{\Emb G} = \frac{\Emb\!_1 \hG}{\Emb \hG}.
$$

Recall that by Corollary \ref{Cr:1Face} we have to count the numberings of
the edges of $G$ that give rise to $1$-faced embeddings. To do this, choose
a cyclic order $\alpha$ corresponding to a $1$-faced embedding.

Let $k_1, \dots, k_m$ be the valencies {\em in the skeleton $\hG$} of the
essential vertices. We can enumerate the numberings of the edges of $G$
that lead to the cyclic order $\alpha$ in the following way. First choose a
set of $k=k_1+ \cdots +k_m$ numbers for the skeleton edges adjacent to the
essential vertices. Second, attribute these numbers to the edges so as to
obtain the desired cyclic order. Third, number the remaining edges with the
remaining numbers (in an arbitrary way). We obtain
$$
\frac{(n+2g-1)!}{k! (n+2g-1-k)!} \cdot \frac{k!}{(k_1-1)! \cdots (k_m-1)!}
\cdot (n+2g-1-k)! = \frac{(n+2g-1)!}{\Emb \hG}.
$$

Thus every cyclic order $\alpha$ contributes by $1/(\Emb \hG)$ to the
coefficient of $G$ in $P_{g,n}/(n+2g-1)!$. Therefore the coefficient of~$G$
equals $\Emb\!_1 \hG/\Emb \hG$. \qed

\bigskip

Now we will evaluate the coefficient of a long graph in~$P_{g,n}$
in another way using Theorem~\ref{Thm:div}.

\begin{definition}
A {\em $k$-spider} $X_k$ is the tree with $k$ vertices
of valency~$1$ and one vertex of valency~$k$.
\end{definition}

\begin{lemma} \label{Lem:onlyspiders}
Disjoint unions of spiders are
the only terms of $R_g$ that give nonzero contributions
to the coefficient of a long graph in the product $T_n R_{g,n}$.
\end{lemma}

\paragraph{Proof.} Let $G$ be a long graph with $n+2g-1$ edges
and $2g$ independent cycles.
On this graph we can mark the $n-1$ edges that come from the
factor $T_n$ (these edges form a tree) and the $2g$ edges
that come from the factor $R_{g,n}$. Since $G$ has $2g$ independent
cycles, erasing each of the $2g$ edges of the second group should
disrupt exactly one cycle of~$G$. It follows that if two edges of
the second group are adjacent, then their common vertex is an
essential vertex of~$G$. Thus the graph formed by the second
group of vertices is composed of isolated edges (which are 1-spiders)
and spiders centered at essential vertices of~$G$. \qed

\begin{remark}
Actually we will soon see that graphs with isolated edges do not
occur in~$R_{g,n}$, so the group of $2g$ edges will never contain
isolated edges.
\end{remark}

 \begin{lemma} \label{Lem:spidercoef}
The coefficient of $X_{2g}$ in $r_g$ and $R_g$ is equal
to $\frac{B_{2g}}{2g}$ and $\frac{1}{2^{2g}(2g+1)}$, respectively.
 \end{lemma}

 \paragraph{Proof.}
According to Definition~\ref{Def:ArR}, we have
\begin{equation} \label{Eq:trace}
r_g = \frac{B_{2g}}{(2g)!(2g)}
\Tr A^{2g} = \frac{B_{2g}}{(2g)!(2g)} \sum_{i_1, \dots, i_{2g}}
A_{i_1 i_2} A_{i_2 i_3} \dots A_{i_{2g-1} i_{2g}} A_{i_{2g} i_1}.
\end{equation}
Obviously, $X_{2g}$ appears only in the terms with $i_1 = \dots = i_{2g}$.
The contribution of these terms is $(2g)!$ because the ``legs''
of the spider can appear in any order. This proves the first statement.

To find the coefficient of $X_{2g}$ in $R_g$ we note that a subgraph of a
spider is always a spider. In other words, the linear span $\mathcal N =
\langle G \mid G \ne X_1, X_2, \dots\rangle \subset \W$ of all the graphs
except spiders is an ideal in $\W$. We can make our computation modulo this
ideal, \ie, neglecting all graphs except the spiders. We have
%*
\begin{equation}\label{Eq:MultN}
\epsfbox[-1 5 6 20]{palochka.eps}^{2g} = (2g)! \; X_{2g} \;\; \bmod
\mathcal N,
\end{equation}
%*
and therefore
 %*
 \begin{equation}\label{Eq:rsmall}
\sum_g r_g = \sum_g \frac{B_{2g}}{2g} X_{2g} =
\sum_g \frac{B_{2g}}{(2g)!(2g)} \;\;
\epsfbox[-1 5 6 20]{palochka.eps}^{2g}
= \ln \phi(\epsfbox[-2 5 5 20]{palochka.eps}) \;\; \bmod \mathcal N,
 \end{equation}
 %*
where, as usual, $\phi(t) = \frac{\sinh (t/2)}{t/2}$. Equations
\eqref{Eq:MultN} and~\eqref{Eq:rsmall} imply
 %*
 \begin{equation*}
\sum_g R_g = \phi(\epsfbox[-2 5 5 20]{palochka.eps})
= \sum_g \frac{1}{2^{2g}(2g+1)!} \;\;
\epsfbox[-1 5 6 20]{palochka.eps}^{2g}
= \sum_g \frac{1}{2^{2g}(2g+1)} X_{2g} \;\; \bmod \mathcal N.
 \end{equation*}
 %*
This proves the second statement. \qed

\begin{lemma} \label{Lem:union}
Let $G$ be a disjoint union of several (more than one)
connected graphs. The coefficient of~$G$ in $r = \sum r_g$
equals~$0$. The coefficient of~$G$ in $R = \sum R_g$
is the product of coefficients of its connected components.
\end{lemma}

\paragraph{Proof.}
The first statement follows from Equation~\eqref{Eq:trace}. Indeed, we have
$(A_n)_{ij} = -w_{ij}$ if $i \ne j$, while $(A_n)_{ii} = \sum_j w_{ij}$.
Thus every term in the sum of Equation~\eqref{Eq:trace} corresponds to a
(posibly self-intersecting) cycle and several edges sticking out of it
(possibly leading to vertices that are already in the cycle). Such a graph
is always connected.

The second statement is an elementary corollary of the first. \qed

\begin{remark}
The last lemma allows us, in particular, to find the coefficient of a
disjoint union of spiders in $R_g$: it is just the product of the
coefficients of the individual spiders. Note that since $r$ contains only
even degree terms, it follows that $R_g$ only contains graphs whose each
connected component has an even degree of edges. In particular, graphs in
$R_g$ do not contain isolated edges. Similarly, the coefficient of the
disjoint union of two $3$-spiders in $R_6$ vanishes.
\end{remark}

\paragraph{Proof of Theorem~\ref{Thm:spiders}.}
Let $G$ be a connected graph with $2g$ independent cycles. Insert a large
enough number of vertices into the edges of~$G$ to obtain a long
graph~$\tG$. This operation changes neither the number of $1$-faced
embeddings, nor the total number of embeddings.

Let $n$ be the number of vertices of $\tG$ and let
us compute the coefficient of $\tG$ in $P_{g,n}/(n+2g-1)!$ using
the right-hand side of the equality of Theorem~\ref{Thm:div}.
We claim that the answer is precisely the sum of weights
of the decorations of~$G$. Indeed, the decorations of~$G$
and the ways to erase $2g$ edges of $\tG$ so as to obtain
a tree are in a natural one-to-one correspondence. Moreover,
the weight of a decoration was defined to coincide with the
coefficient of the corresponding disjoint union of spiders
in $R_g$.

On the other hand, the coefficient of $\tG$ in $P_{g,n}/(n+2g-1)!$
equals $\Emb\!_1 \tG/ \Emb \tG$ according to
Proposition~\ref{Prop:embeddings}.

Thus the sum of weights of the decorations of~$G$ equals
$$
\frac{\Emb\!_1 \tG}{ \Emb \tG} = \frac{\Emb\!_1 G}{ \Emb G}.
$$
\qed

\section*{Appendix: Hurwitz generating function and the cut-and-join
equation}

In this Appendix we consider coverings that are not pseudo-polynomial (but
are still Morse except at infinity). Namely, fix a partition $\lambda =
(\lambda_1, \dots, \lambda_s)$ (where $\lambda_1 \le \dots \le \lambda_s$)
of the number 
%$n \bydef \lmod \lambda\rmod \bydef \lambda_1 + \dots + \lambda_s$ 
$n = \lambda_1 + \dots + \lambda_s$ 
into $s$
%$s \bydef \#\lambda$ 
parts; consider degree~$n$ (not
necessarily connected) ramified coverings of the sphere by a surface $M_g$
of Euler characteristics $\chi(M_g) = 2-2g$ with simple ramifications over
$n+2g-2+s$ points and a ramification of type $\lambda$ over one point
labeled $\infty$; the latter condition means that $\infty$ has $s$
preimages with multiplicities $\lambda_1, \dots, \lambda_s$. Since $M_g$ is
not assumed connected, it is possible to have $g < 0$.

The monodromy of the covering in question over~$\infty$ is an element of
the conjugacy class $C_\lambda \subset S_n$ that consists of permutations
that are products of $s$ independent cycles of lengths $\lambda_1, \dots,
\lambda_s$. The monodromies over the other branch points are
transpositions. Similarly to Definitions \ref{Df:CycleFact} and
\ref{Df:HurwPoly} a list of transpositions $\tau_1, \dots, \tau_{n+2g-2+s}
\in S_n$ is called a {\em genus $g$ factorizion of $\lambda$} if
$\tau_{n+2g-2+s} \cdots \tau_1 \in C_\lambda$. Consider the {\em Hurwitz
polynomial}
 %*
 \begin{equation*}
P_{g,\lambda} \bydef \sum w(\tau_1) \cdots w(\tau_{n+2g-2+s}),
 \end{equation*}
 %*
where the sum is taken over all such factorizations; recall that the
numbers $n$ and $s$ are determined by~$\lambda$.

The {\em Hurwitz number} $h_{g,\lambda}$ is defined by $h_{g,n} =
\frac1{n!} P_{g,\lambda}(\bf 1)$.

The condition $\tau_{n+2g-2+s} \cdots \tau_1 \in C_\lambda$ is invariant
under conjugation in $S_n$, so that $P_{g,\lambda}$ is an $S_n$-invariant
polynomial. Its degree is $m \bydef n+2g-2+s$, so, once $\lambda$ and
$P_{g,\lambda}$ are known, it is possible to determine $g$. Therefore one
can collect all the Hurwitz polynomials into a general {\em Hurwitz
$\w$-generating function}
 %*
 \begin{equation}\label{Eq:HurwW}
H(\w; p_1, p_2, \dots) = \sum_{g,\lambda} \frac{1}{m!}
\frac{P_{g,\lambda}(\w)}{n!} \frac{p_{\lambda_1} \dots p_{\lambda_s}
}{|\name{Aut}(\lambda)|}\in \C[[\w, p_1, p_2, \dots]];
 \end{equation}
 %*
here $\w$ means the infinite collection of variables $w_{ij}$, $1 \le i <
j$, and $\name{Aut}(\lambda)$ is the set of permutations $\sigma$ of
$s$ elements such that $\lambda_i = \lambda_{\sigma(i)}$ for all $i$. Taking
$w_{ij} = \beta$ for all $i,j$ one obtains the usual Hurwitz generating
function
 %*
 \begin{equation}\label{Eq:HurwNormal}
h(p_1, p_2, \dots) = \sum_{g,\lambda} \frac{\beta^m}{m!} h_{g,\lambda}
\frac{p_{\lambda_1} \dots p_{\lambda_s}}{|\name{Aut}(\lambda)|}
 \end{equation}
 %*
studied intensively during the last decade (see e.g.\
\cite{GJ,GJV,WitSm,PZv}, among many others). One of the main properties of
the Hurwitz generating function \eqref{Eq:HurwNormal} is that it satisfies
the so-called cut-and-join equation \cite{GJ}; here we generalize it to the
$\w$-function \eqref{Eq:HurwW}.

Define the {\em cut-and-join operator} by
 %*
 \begin{equation*}
L = \frac{1}{2} \sum_{k,l \geq 1} \bigl((k+l) p_k p_l \pder{}{p_{k+l}} +
kl p_{k+l} \pdertwo{}{p_k}{p_l}\bigr)
 \end{equation*}
 %*

 \begin{theorem}
The Hurwitz generating $\w$-function $H$ satisfies the {\em cut-and-join
differential equation}
 %*
 \begin{equation}\label{Eq:CutAndJoin}
\sum_{1 \le i < j} \pder{H}{w_{ij}} = LH.
 \end{equation}
 %*
 \end{theorem}

\paragraph{Proof.} Let $G$ be a graph without loops with $n$ numbered
vertices and $m$ edges. Let $\lambda$ be
a Young diagram with $n$ squares.

Let $S_1(G,\lambda)$ be the set of pairs $(e,u)$ where
$e$ is an additional edge joining two vertices of~$G$ and
$u$ is a numbering of the set $\{ \mbox{edges of}~$G$ \} \cup
\{e\}$ such that the product $\rho$ of transpositions corresponding to the
numbering lies in the conjugacy class~$C_\lambda$.

Futher, let $S_2(G,\lambda)$ be the set of pairs $(v, \{i,j\})$ where $v$ 
is a numbering of the edges of~$G$ and $\{i,j\} \subset \{1, \dots n\}$ is 
an unordered pair of elements such that $\sigma \tau_{ij}$ lies in the 
conjugacy class~$C_\lambda$ (where $\sigma$ is the product of 
transpositions corresponding to the numbering).

The coefficient of the monomial $w_G p_\lambda$
(with self-explanatory notation) in the left-hand side
of~\eqref{Eq:CutAndJoin} equals $\lmod S_1(G,\lambda)\rmod/n!(m+1)!
\lmod\Aut(\lambda)\rmod$. The coefficient of $w_G p_\lambda$ in the
right-hand side equals $\lmod S_2(G,\lambda)\rmod/n!m!
\lmod\Aut(\lambda)\rmod$. (More
precisely, if $i$ and $j$ lie in two different cycles of $\sigma$ of
lengths $k$ and $l$, this corresponds to the term $kl p_{k+l}
\partial^2/\partial p_k \partial p_l$; if $i$ and $j$ lie in the same cycle
of $\sigma$ of length $k+l$, $k$ and $l$ being the distances from $i$ to
$j$ and from $j$ to $i$ along the cycle, this corresponds to the term
$(k+l) p_k p_l \partial/\partial p_{k+l}$.)

It remains to establish an $(m+1)$-to-$1$ correspondence between $S_1$ and
$S_2$. For a pair $(e,u) \in S_1$ we
write $\rho = \rho_2 \tau(e) \rho_1$. Here $\rho_1$ is the product of the
tranpositions coming before $e$, $\tau(e)$ is the transposition
corresponding to $e$, and $\rho_2$ is the product of transpositions coming
after $e$. Then we let $\sigma = \rho_2\rho_1$, which also determines a
numbering $v$ of the edges of~$G$; $\{i,j\}$ are the elements permuted by
$\tau(e)$. The permutation $\sigma \tau_{ij} = \rho_2\rho_1\tau(e)$
is conjugate to $\rho$ and therefore lies in the
class~$C_\lambda$. Finally, the map from $S_1$ to
$S_2$ we have constructed is invariant under circular permutations of the
numbering~$u$, thus it is actually an $(m+1)$-to-$1$ correspondence. \qed

\begin{remark}
Note that the last part of the proof is different from the usual proof of
the classical cut-and-join equation, where one can always assume that
$\tau(e)$ is the first transposition in the list and there is no need to
introduce the cirular change in the numbering of transpositions.
\end{remark}


\begin{thebibliography}{99}

\bibitem{BouMar} {\bf V.\ Bouchard, M.\ Mari\~no.} {\em Hurwitz numbers,
matrix models and enumerative geometry.} -- {\tt arXiv: 0709.1458}
(21~pp.).

\bibitem{ELSV} {\bf T.\ Ekedahl, S.\ Lando, M.\ Shapiro, A.\ Vainshtein.}
{\em Hurwitz numbers and intersections on moduli spaces of curves.}
-- Invent.\ Math.~{\bf 146} (2001),
no~2, pp.~463--533. {\tt arXiv: math/0004096}.

\bibitem{FuHa} {\bf W.\ Fulton, J.\ Harris.} {\em Representation theory. A
first course.} -- Springer, Berlin, 1991.

\bibitem{GJ} {\bf I.P.\ Goulden, D.M.\ Jackson.} {\em Transitive
factorisation into transpositions and holomorphic mappings on the sphere}
-- Proc.\ Amer.\ Math.\ Soc., {\bf 125} (1997), no.~1, pp.~51--60.

\bibitem{GJV} {\bf I.P.\ Goulden, D.M.\ Jackson, R.\ Vakil.} {\em The
Gromov--Witten potential of a point, Hurwitz numbers, and Hodge integrals}
-- Proc.\ London Math.\ Soc., {\bf 83} (2001), no.~3, pp.~563--581.

\bibitem{GouldYong}{\bf I.\ Goulden, A.\ Yong.} {\em Tree-like properties
of cycle factorizations.} -- J.\ Comb.\ Theory, Ser.~A, {\bf 98} (2002),
no.~1, pp.~106--117.

\bibitem{GroTay} {\bf D.\ Gross, W.\ Taylor.} {\em Two-dimensional QCD is a
string theory.} -- Nucl.\ Phys.~B~{\bf 400} (1993), no~1--3, pp.~181--208.

\bibitem{HarZag} {\bf J.\ Harer, D.\ Zagier.} {\em The Euler characteristic
of the moduli space of curves.} -- Invent.\ Math.~{\bf 85} (1986), no.~3,
pp.~457--485.

\bibitem{WitSm} {\bf M.E.\ Kazarian, S.K.\ Lando.} {\em An
algebro-geometric proof of Witten's conjecture} -- J.\ Am.\ Math.\ Soc.,
{\bf 20} (2007), no.~4, pp.~1079--1089.

\bibitem{Kirchhoff} {\bf G.\ Kirchhoff.} {\em \"Uber die Aufl\"osung der
Gleichungen, auf welche man bei der Untersuchung det linearen Verteilung
galvanischer Str\"ome gefurht wird} -- Ann.\ Phys.\ Chem., {\bf 72} (1847),
pp.~497--508.

\bibitem{LanZvo} {\bf S.K.\ Lando, A.\ Zvonkin.} {\em Graphs on Surfaces
and their applications.} With an appendix by D.~Zagier. Encyclopaedia of
Mathematical Sciencies~{\bf 141}. Low-dimensional topology,~II.
Springer-Verlag, Berlin, 2004.

\bibitem{Macdonald} {\bf I.G.\ Macdonald.} {\em Symmetric functions and
Hall polynomials} -- Oxford Mathematical Monographs, Clarendon Press,
Oxford, 1979.

\bibitem{PZv}  {\bf D.\ Panov, D.\ Zvonkine.} {\em Enumeration of almost
polynomial rational functions with given critical values} -- Eur.\ J.\ Comb.\
{\bf 29} (2008), no.~2, pp.~470--479.

\bibitem{PouSch} {\bf D.\ Poulalhon, G.\ Schaeffer.} {\em Factorizations of
large cycles in the symmetric group.} -- Discrete Math.~{\bf 254} (2002),
no.~1--3, pp.~433--458.

\bibitem{ShShVa} {\bf B.\ Shapiro, M.\ Shapiro, A.\ Vainshtein.} {\em
Ramified coverings of of~$S^2$ with one degenerate branching point and
enumeration of edge-ordered graphs.} -- Amer.\ Math.\ Soc.\ Transl.\
Ser.~2.,~{\bf 180} (1997), Topics in Singularity Theory, no.~42,
pp.~219--227. {\tt www2.math.su.se/\~{}shapiro/Articles}
 \end{thebibliography}
\end{document}